# On Self-Descriptive Squares

Lee Sallows and Dmitry Kamenetsky

*An innocent seeming curio develops some surprisingly complex ramifications*

**Abstract:** A novel kind of self-referential square matrix is introduced. A certain subset of the matrix entries record the frequencies of occurrence of each distinct number appearing within the entire matrix. Such squares are necessarily elusive. Our investigation brings to light interesting cases, such as 'generic' squares that include algebraic variables and self-descriptive squares that are also magic squares.

We call a square array of numbers, usually integers, "self-descriptive" (s-d, for short) when the sum of the entries in any row or column describes the number of times that a row's rightmost entry, or a column's lowermost entry, appear in the square. Figure 1 shows an example of size 4×4, where the sums of the rows and columns are shown at right and below. Here the sum of the numbers in the top row is 3, which is the number of times that the integer 2 occurs in the square, while the sum of the numbers in the righthand column is 1, being the number of times that 4 appears. And so on for the remaining rows and columns. We call the number of occurrences of a number in an array its *frequency*.

| 2 | 1 | -2 | 2 | 3 |
|---|---|---|---|---|
| 1 | 2 | 3 | -4 | 2 |
| 3 | 1 | -1 | -1 | 2 |
| -4 | -2 | 3 | 4 | 1 |
| 2 | 2 | 3 | 1 | |

Figure 1. A trivial 4×4 s-d square.

Two distinct flaws mar the elegance of Figure 1. Firstly, no row or column sum identifies the frequency of integer 1, which is nevertheless present in the square. Secondly, the frequency of -4 is needlessly identified twice: once in the second row and once in the first column. The self-description is thus both incomplete and redundant.

It is for this reason that a self-descriptive square of order *n* (or size *n*×*n*) is deemed *trivial* unless every distinct integer appearing in the square occurs exactly once among the $2n-1$ cells making up the shaded L-shaped border (hereinafter simply 'the border') formed by the right-hand column and bottom row. In this way, with the single exception of the twice-counted lower right-hand corner entry (see box on page 2), the frequency of every distinct integer appearing in the square is identified by a unique row or column sum. Figure 2 shows an example of a non-trivial self-descriptive square, the number of distinct integers occurring now being $2n-1 = 7$.

| 1 | 3 | 3 | -4 | 3 |
|---|---|---|---|---|
| 3 | -4 | 2 | 1 | 2 |
| -1 | 2 | -4 | 4 | 1 |
| -1 | 0 | 2 | 3 | 4 |
| 2 | 1 | 3 | 4 | |

Figure 2. A non-trivial s-d square.

That the frequency of the lower right-hand corner number is given both by the bottom row and right-hand column sums is regrettable, because it means that the 2n row/column sums identify the frequencies of only 2n-1 distinct numbers rather than 2n. However, this is the price paid for a avoiding a different kind of blemish. Take for example, the square shown in the Figure below that is based upon a quite different self-descriptive scheme. Here the r-th row sum describes the frequency of its r-th entry (from left to right), while the c-th column sum gives the frequency of its c-th entry (from top to bottom). The enumerated entries thus become those lying along the two main diagonals, shown shaded. For example, the sum of the third row is 4, which is the frequency of -1, its third entry; while the sum of the second column is 3, which is the frequency of -3, its second entry.

| -4 | 4  | -3 | 4  | 1 |
|----|----|----|----|---|
| 4  | -3 | 2  | -1 | 2 |
| 3  | 3  | -1 | -1 | 4 |
| -2 | -1 | 3  | 1  | 1 |
| 1  | 3  | 1  | 3  |   |

The 2n row/column sums now describe the frequencies of 2n distinct numbers. Problem solved? For squares for which n is even, yes. But when n is odd we find that the two diagonals intersect on the centre cell, giving us a total of only 2n-1 entries. Similar shortcomings arise with alternative schemes. An advantage of the definition settled on here is that the number enumerated by each row or column sum is unambiguously located at a fixed position at the end of that row or column without need of further consideration. It is the lack of such simplicity and clarity that is the great drawback of alternative definitions.

Optional box referred to in previous paragraph.

The solution represented by Figure 2 is far from unique. A moment's thought shows that its three upper rows could be reshuffled in any order without impairing the square's self-descriptive property, and likewise its three left-hand columns. The same applies more generally to squares of size $n \times n$: excepting the row and column that together form the border, the remaining $n$-1 rows and/or $n$-1 columns can be permuted in any order with impunity. Moreover, a square remains self-descriptive when reflected about its main (\) diagonal, a change that merely switches the entries in the $m^{th}$ row with those in the $m^{th}$ column. The total number of distinct possible rearrangements is therefore: number of row permutations × number of column permutations × number of reflections = $(n-1)!^2 \times 2$. Figure 3 shows one of the $3!^2 \times 2 = 72$ self-descriptive rearrangements of Figure 2, one that results in the row and column sums (which happen to be the same in this case) now appearing in order of magnitude. We say that two squares that are derivable from eachother by means of these s-d-preserving transformations are *equivalent*.

| 2 | -1 | -4 | 4 | 1 |
|---|----|----|---|---|
| -4 | 3 | 2 | 1 | 2 |
| 3 | 1 | 3 | -4 | 3 |
| 0 | -1 | 2 | 3 | 4 |
| 1 | 2 | 3 | 4 | |

Figure 3.
A rearrangement of the square from Figure 2.

It is natural to think of Figures 2 and 3 as being alternative presentations of a basic underlying square. The latter can only be one of the 72 possibilities, but how do we decide which? One answer is to identify it as that (unique) specimen which is expressed in *standard normal form* (SNF). Such a square can be defined as one in which:
- The first $n$-1 columns are sorted such that the first $n$-1 entries in row $n$ appear in increasing order of magnitude from left to right
- The first $n$-1 rows are sorted such that the first $n$-1 entries in column $n$ are sorted in increasing order of magnitude from above to below
- The top-right entry > bottom left entry (since both entries are in the border they must be distinct)

Returning now to Figures 2 and 3, we find the latter meets only the third requirement above, while Figure 2 satisfies both of the first two requirements, but not the third. On reflecting Figure 2 about its main (\) diagonal however, we correct this shortcoming to arrive at Figure 4, which is thus the only square among the 72 to satisfy the demands of standard normal form. Similarly, the rows and columns of Figure 3 can also be permuted to result in Figure 4.

| 1 | 3 | -1 | -1 | 3 |
|---|---|----|----|---|
| 3 | -4 | 2 | 0 | 2 |
| 3 | 2 | -4 | 2 | 1 |
| -4 | 1 | 4 | 3 | 4 |
| 2 | 1 | 3 | 4 | |

Figure 4.
A square in standard normal form.

We call Figures 2 and 3 *anagrams* of eachother because they are composed of the same set of $n^2$ entries, with each entry occurring the same number of times. The authors suspected that any such anagrammatic pair must therefore be equivalent, an assumption that turned out to be incorrect. At left in Figure 5 we see a pair of 5×5 s-d squares composed of the same set of 25 entries. On the right, each square is re-expressed in SNF, but differ from eachother. The two squares at left are thus non-equivalent after all.

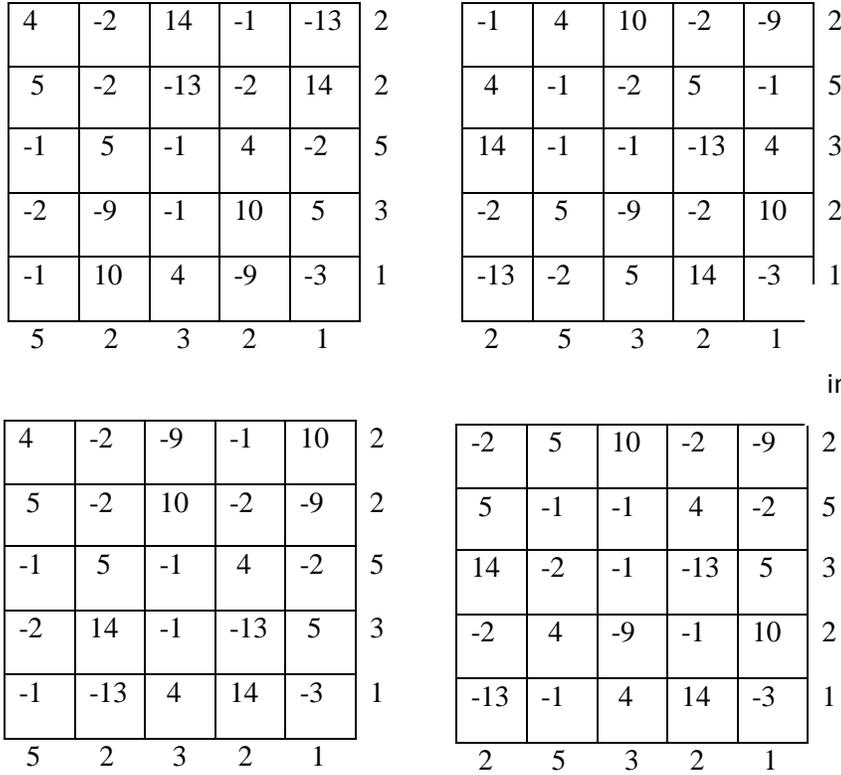

Figure 5.
Re-expressing two anagrams in SNF reveals their non-equivalence.

*Self-descriptive squares of low order*

There exists just one s-d square of size 1×1, which consists of a single cell containing the number 1.

Since, in any s-d square, the bottom row and right-hand column sums are equal, in a 2×2 square this will mean that the top-right and bottom-left entries must be the same. But both occupy positions in the border, and as such must be distinct if the square is to be non-trivial. It is easy to see that, even then, the only integers able to satisfy the self-descriptive requirement are as shown in Figure 6.

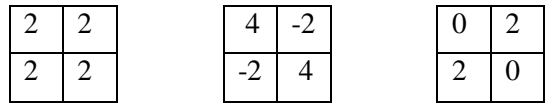

Figure 6.
Three trivial 2×2 squares.

Clearly the row and column sums in any s-d square must all be integers $\geq 1$, but we can say a little more. Let $\Sigma r_m$ and $\Sigma c_m$ stand for the sum of the entries in the $m^{th}$ row and $m^{th}$ column, respectively. Then for any s-d square $S$ of size $n \times n$, we have: $\Sigma r_1 + \Sigma r_2 + .. + \Sigma r_n$ = sum of the $n^2$ entries in $S$ = $\Sigma c_1 + \Sigma c_2 + .. + \Sigma c_n$.

Further, as already noted, $\Sigma r_n = \Sigma c_n$, from which $\Sigma r_1 + \Sigma r_2 + .. + \Sigma r_{n-1} = \Sigma c_1 + \Sigma c_2 + .. + \Sigma c_{n-1}$. However, since $\Sigma r_1 + \Sigma r_2 + .. + \Sigma r_n + \Sigma c_1 + \Sigma c_2 + .. + \Sigma c_n = n^2 + \Sigma r_n$ (the $n^2$ entries of $S$, but with $\Sigma r_n = \Sigma c_n$ counted twice), on subtracting $\Sigma r_n$ and $\Sigma c_n$ from both sides, we have: $\Sigma r_1 + \Sigma r_2 + .. + \Sigma r_{n-1} + \Sigma c_1 + \Sigma c_2 + .. + \Sigma c_{n-1} = n^2 - \Sigma r_n$, or

$$2 \times (\Sigma r_1 + \Sigma r_2 + .. + \Sigma r_{n-1}) = 2 \times (\Sigma c_1 + \Sigma c_2 + .. + \Sigma c_{n-1}) = n^2 - \Sigma r_n \qquad \text{[Equ. 1]}$$

from which is seen that $n^2 - \Sigma r_n$ is divisible by 2, implying that $\Sigma r_n = \Sigma c_n$ and $n$ must be of same parity. Hence in any s-d square of size $n \times n$, the bottom row sum (= right-hand column sum) will be even only when $n$ is even, and odd only when $n$ is odd. From this is seen that in any s-d square of $3 \times 3$, $\Sigma r_3 = \Sigma c_3$ is an odd integer $\geq 1$. However, we can go a little further and show that $\Sigma r_3 = \Sigma c_3 \leq 3$.

Suppose $\Sigma r_3 = \Sigma c_3 = 5$, the least possible odd integer $>3$. Then from [Equ. 1] above, $\Sigma r_1 + \Sigma r_2 = \Sigma c_1 + \Sigma c_2 = \frac{1}{2}(n^2 - \Sigma r_3) = 2$. Hence, since they must both be integers, $\Sigma r_1 = 1$ and $\Sigma r_2 = 1$, as shown in Figure 7 in which $v$, $w$, $x$, $y$ and $z$ represent the five distinct border entries.

|   |   | $v$ | 1 |
|---|---|---|---|
|   |   | $w$ | 1 |
| $z$ | $y$ | $x$ | 5 |
| 1 | 1 | 5 |   |

Figure 7.

However, since the frequency of $x$ is 5, the 4 empty cells in the 2×2 subsquare must contain the 4 missing $x$'s. But then, from the two top rows we would have $2x + v = 1$ and $2x + w = 1$, implying $v = w$, contrary to our demand they be distinct. That the frequency of $x$ could not be 5 or greater is now obvious. Thus $\Sigma r_3 = \Sigma c_3$ can only be 1 or 3.

Continuing in similar vein, we leave it as an exercise for the reader to complete the (admittedly trivial, yet time and space-consuming) proof that the two squares of Figure 8 turn out to be the only existent s-d squares of order 3, shown in SNF.

| 2 | -2 | 2 | 2 |
|---|---|---|---|
| 0 | -2 | 3 | 1 |
| 0 | 5 | -2 | 3 |
| 2 | 1 | 3 |   |

| 2 | 0 | 0 | 2 |
|---|---|---|---|
| 2 | -2 | 1 | 1 |
| -2 | 3 | 2 | 3 |
| 2 | 1 | 3 |   |

Figure 8.
The only two 3×3 s-d squares.

We have already glanced at a couple of s-d squares of orders 4 and 5. In the absence of deeper insight, at present the only known way of producing squares larger than $3 \times 3$ is via computer searches. In this way we are able to report that, for example, there exist at least 28 distinct squares of order 4 using integers with absolute value $\leq 4$. Similarly, there exist at least 720 squares of order 5 using integers with absolute value $\leq 4$, and at least 33305 squares of order 6 using integers with absolute value $\leq 5$. Such data is of course of limited value, and particularly so in light of the unexpected discovery of *generic* s-d squares (for which see below) bringing to light the existence of families of squares that are infinite in number.

*Computer searches for s-d squares*

The interwoven reflexive properties of s-d squares render them peculiarly elusive objects.  We used a computer program to find most of our examples.  The algorithm is based on the hill-climbing method, the aim being to incrementally improve the quality of a partial solution until a full solution is found [1][2]. Our approach is outlined in the pseudocode shown in Figures 9 and 10.

```
int[][] solve(int N, int M)
{
   numberOfMutations = 3    //mutations performed in each iteration
   bestErrors = Infinity   //least number of errors seen so far

   //randomly initialize an NxN square with values in [-M,+M]
   bestSquare = initialize(N, M)
   while True:
   {
     //mutate the current best square
     square = mutate(bestSquare, numberOfMutations)

     //optimize the mutated square
     optimize(square)
     e = countErrors(square)
     if e < bestErrors
     {
       bestErrors = e
       bestSquare = square
       if bestErrors == 0 then return bestSquare   //solution found!
     }
   }
}
```

Figure 9. Main routine of the algorithm.

An s-d square must comply with the following four conditions:
1. Every border number must appear exactly once in the border.
2. Every number appearing in a non-border cell must also appear in a border cell.
3. Every row sum must be equal to the number of times its rightmost number occurs in the square.
4. Every column sum must be equal to the number of times its bottom number occurs in the square.

We measure the quality of a partial solution by counting the number of times that any one of these four conditions is infringed (this is performed in the *countErrors* routine). As the number of violations (errors) decreases, so the quality of a solution increases. A full solution is obtained once the number of violations has dropped to zero.

The inputs to the algorithm (see the *solve* routine in Figure 9) are the size of the square (N) and the maximum absolute value of its numbers (M). The algorithm begins by creating a candidate square (approximation to a solution) with randomly assigned values in the range [-M, +M] (performed in the initialize routine). An iterative process is now started in which the algorithm performs some number of mutations (we use 3) to the current best candidate square. A mutation involves selecting a random cell and replacing it with a new value chosen at random in the range [-M, +M]. An attempt is made to improve the quality of the mutated solution via the *optimize* routine shown in Figure 10. This routine makes small changes in an attempt to improve the quality of the current solution; if a change is beneficial then it is accepted, otherwise it is rejected. As with the mutations, these adjustments involve changing a single cell to a new value. If the optimized square is of higher quality than the best square found so far, then it becomes the best candidate square for use in subsequent iterations. Finally the algorithm terminates when the number of errors in a solution reaches 0.

```
void optimize(int[][] square, int N, int M)
{
  bestErrors = countErrors(square)
  for every cell (r,c) and value m in [-M,+M]:
  {
    oldValue = square[r][c]   //save the old value
    square[r][c] = m  //make a temporary assignment

    e = countErrors(square)
    //keep assignment if it reduces errors, otherwise undo
    if e < bestErrors
      bestErrors = e
    else
      square[r][c] = oldValue
  }
}
```

Figure 10. The optimize routine for improving the quality of partial solutions.

Al Zimmermann pointed out to us that the *optimize* routine can be made substantially faster. Instead of recalculating the number of errors for each change in value (using *countErrors*) in $O(n^2)$ time one can do it in $O(1)$ time. Each *tuple* in the square (that is, each row and each column) makes an additive contribution to the number of errors. When a single cell is modified, 4 tuples change the amount of their contribution. For example, if the value at column c, row r is changed from $v_0$ to $v_1$, then the tuples whose contributions change are column c, row r, and the two tuples whose last cells have values $v_0$ and $v_1$. To calculate the change in value for the objective function, we need only calculate the changes in contribution from these 4 tuples. It must be noted that the 4 tuples are not necessarily distinct and therefore the number of tuples to be handled could be fewer than 4. In order to calculate the modified number of errors we must know: the current number of errors, the sum of values for each tuple of the square, and the count of occurrences for each value in the square. Some special handling is required if $v_0$ or $v_1$ is the bottom-right number in the square.

*Generic s-d squares*

We call an s-d square *generic* when it contains one or more algebraic variables for which numerical values can be substituted. Hence a generic square represents a potentially infinite family of distinct numerical s-d squares. There exist no generic squares of order 3, but Figure 11 shows six examples of order 4. Our computer investigations suggest that they can be found for all higher orders.

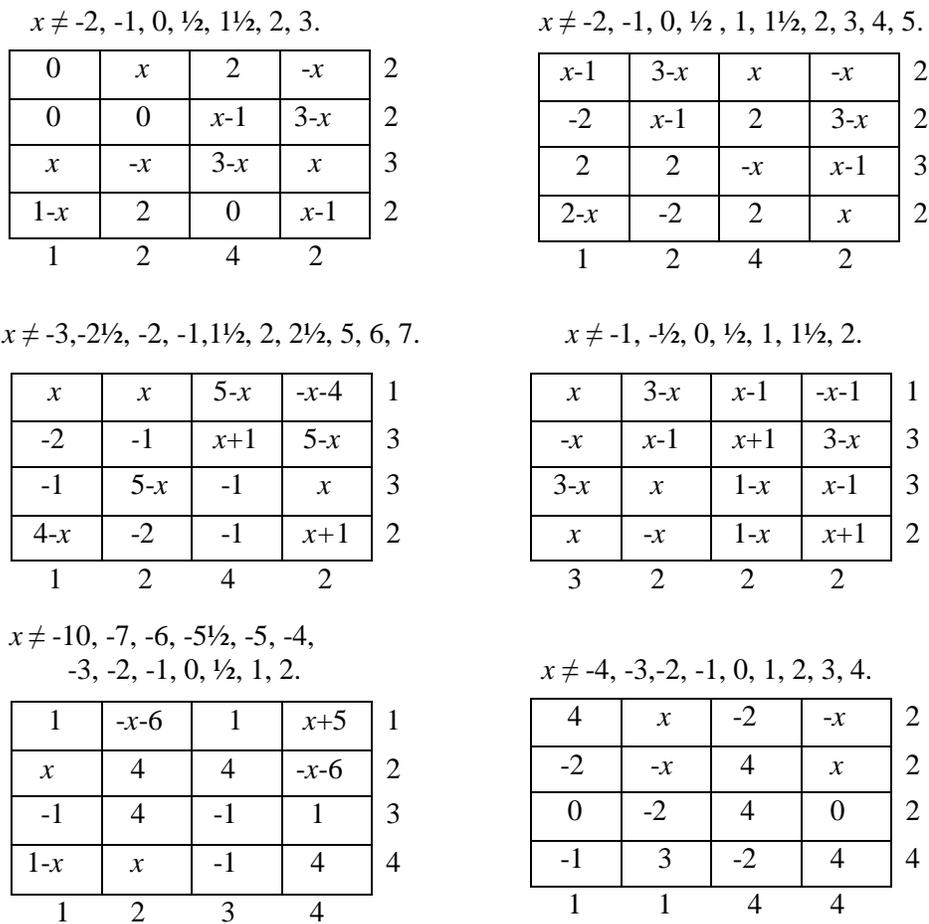

Figure 11.
Six 4×4 generic squares.

A feature common to all generics is that for every variable $v$ occurring in row $r$ and column $c$, the term $-v$ will be found to occur at least twice, once elsewhere in row $r$ and once elsewhere in column $c$. In this way, $v$ and $-v$ always cancel eachother, their net contribution to their row and column sums thus being zero. So, since $-(-v) = v$, a variable and its negative must each appear at least twice in the square. Among the 4×4 examples of Figure 11, the single variable $x$ appears either 2, 3, 5 or 8 times.

A surprising consequence of this cancellation of values is that a number substituted for $v$ need not be integral, or even real; it could even be complex.

In light of these remarks it may be thought that the number chosen to replace a variable in a generic square could be of *any* value. This however, is not the case; certain values are *forbidden*. Suppose, for example, that in the top left square of Figure 11 the variable $x$ is assigned the value -2. But then the value of the term $-x$, which appears twice in the square, will become 2, meaning that 2 will now occur four times in the square instead of twice, as represented by the sum of column 2. Forbidden values may not always be whole numbers. In the square just considered, setting $x$ equal to ½ makes the entry 1-$x$ also equal to ½, thereby rendering the sum of row 3 a false total. The forbidden values for each of the squares of Figure 11 are listed above each square. In general, variables cannot take those values that would make any two border cells equal. Hence the maximum number of forbidden values that can be associated with a single variable is the number of distinct pairs of border cells, given by $(2n-1)(2n-2)/2$.

Incidentally, there is even a trivial generic square of order 2 for which the forbidden value is 1; see Figure 12.

$x \neq 1$

|   |     |   |
|---|-----|---|
| x | 2-x | 2 |
| 2-x | x | 2 |
| 2 | 2 |   |

Figure 12.
A trivial square of order 2.

The generic squares looked at above have all contained a single independent variable. Thus far we have discovered a few generics of orders 5 and 6 using two variables, but none as yet with more; see Figure 13. Note that Figure 5 top left is a numerical instance of the 5×5 square of Figure 13 in which $x = -13$.

| 4 | -2 | 1-x | -1 | x | 2 |
|---|----|-----|----|----|---|
| 5 | -2 | x | -2 | 1-x | 2 |
| -1 | 5 | -1 | 4 | -2 | 5 |
| -2 | y-6 | -1 | 7-y | 5 | 3 |
| -1 | 7-y | 4 | y-6 | -3 | 1 |
| 5 | 2 | 3 | 2 | 1 |   |

| 3+x | 4+x | 2 | -2-x | -5-x | 2 | 4 |
|-----|-----|---|------|------|---|---|
| -2-x | 4+x | 4+x | 6-y | -3-x | y-7 | 2 |
| 2 | -2-x | 4+x | 1 | -3-x | 4+x | 6 |
| 1 | 1 | -2-x | y-7 | 3+x | 6-y | 2 |
| 4+x | 0 | -5-x | 2 | 5+x | -3-x | 3 |
| -5-x | -2-x | 1 | 3+x | 5+x | 0 | 2 |
| 3 | 5 | 4 | 3 | 2 | 2 |   |

Figure 13. Two examples of generic squares showing two distinct variables.

Whether or not there exist higher order generics exhibiting still greater numbers of independent variables is unknown. Incidentally, it should be clear that generic squares elude expression in SNF.

Generic s-d squares can sometimes be discovered through modifying a numerical s-d square. Consider for example the s-d square at left in Figure 14. On replacing every 5 with $x$ and every -5 with $-x$ we obtain the generic square at right. The transformation works because every row or column containing a 5 is accompanied by a single -5 in the same row or column, and thus cancel eachother out.

| 4 | 5 | -2 | -5 | 2 |
|---|---|----|----|---|
| -2 | -5 | 4 | 5 | 2 |
| 0 | -2 | 4 | 0 | 2 |
| -1 | 3 | -2 | 4 | 4 |
| 1 | 1 | 4 | 4 |   |

| 4 | x | -2 | -x | 2 |
|---|---|----|----|---|
| -2 | -x | 4 | x | 2 |
| 0 | -2 | 4 | 0 | 2 |
| -1 | 3 | -2 | 4 | 4 |
| 1 | 1 | 4 | 4 |   |

Figure 14. Deriving a generic square from a numerical specimen.

More complex transformations can also yield generic squares. Consider for example Figure 15(a). Replacing 4 with $x$ and -4 with $-x$ leads to Figure 15(b), a square that is not yet fully generic because its third row and first column show a surplus $x$ in their sums. This can be rectified by replacing further numbers with variables, the aim being to make changes that will cancel out the $x$'s in every row and column. One possible change is to replace 3 with $x$-1 (noting that 3 = 4-1), -3 with 1-$x$ ( -3 = 1-4) and -1 with 3-$x$ (-1 = 3-4) to obtain the generic square shown in Figure 15(c). Moreover, we can generate a distinct generic square simply by choosing a different initial replacement. For example, on changing each 3 to *x* (so that 4 becomes *x*+1) we obtain the generic square of Figure 15(d).

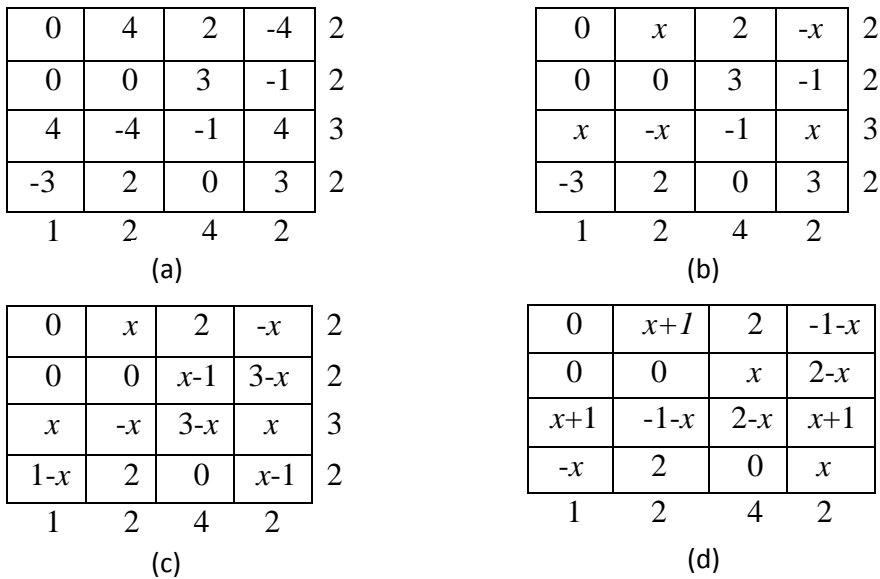

Figure 15. More complex exchanges again yield generic squares.

The algorithm we use for finding such generic squares is as follows. An s-d square is generated using the routine *solve,* described in the *Computer searches* section above. All possible contingent variable replacements that might result in a generic square are then tested in turn. If no successful replacements are found, a new s-d square is then generated for testing. This process is repeated until a valid generic square is found. For a square of order *n*, there are $2n$-1 unique values that can be replaced. Each value can either be replaced or not, resulting in a total of $2^{2n-1}$ cases to be checked.

*Exotic s-d squares*

Given the interlocking requirements that must be satisfied by the numbers occurring in any self-descriptive square, our discovery of more exotic specimens exhibiting still further reflexive properties came as a surprise. In this respect, s-d squares share something in common with magic squares, a similarity that suggested the term *magic* for s-d squares such as those in Figure 16 in which, besides rows and columns, the two main diagonals are also self-descriptive. In this example it is the bottom left and bottom right corner numbers whose frequencies are identified by the main (\) and co-diagonal (/) sums, but squares using alternative corner pairs are also to be found.

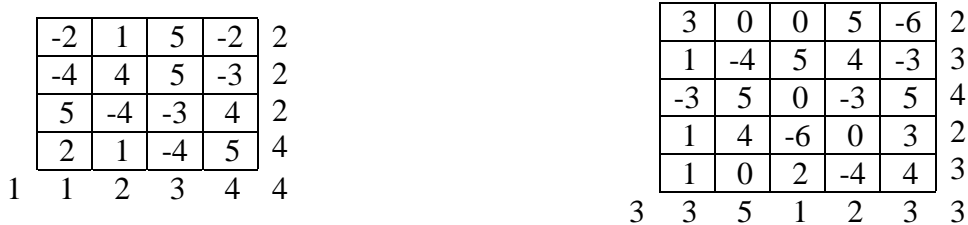

Figure 16. Two examples of magic s-d squares.

We call a row or column *bi-directional* when the frequency of its first entry is the same as that of its last. The frequency of either element is thus given by the sum of the entries in that row or column. This will obviously be the case if the two entries are identical, but may also be true when they are distinct. Figure 17 presents two examples of s-d squares in which every row and column is bi-directional, as indicated by the repeated row and column totals shown.

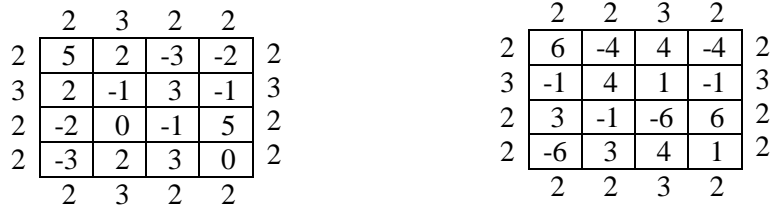

Figure 17. Two s-d squares with bi-directional rows and columns.

Likewise, self-descriptive diagonals may also be bi-directional. In which case, the question of *which* of its two corner elements is to be referred to by the diagonal sum, disappears. This is to be seen in what we call *perfect* s-d squares, being magic s-d squares showing bi-directional rows, columns *and* diagonals. Figure 18 shows two examples.

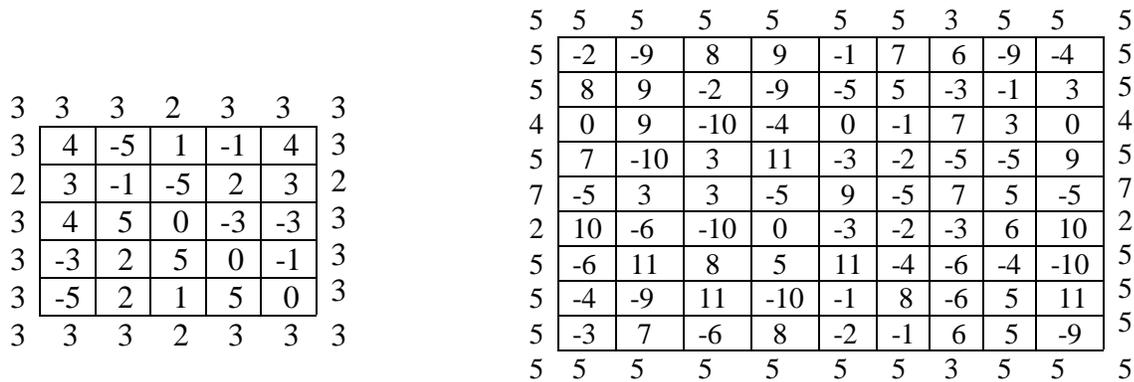

Figure 18. Examples of perfect s-d squares: magic, plus bi-directional orthogonals and diagonals. Left: smallest known perfect s-d square. Right: largest known perfect s-d square.

Also worth mentioning here are *minimal* s-d squares of *n×n* whose entries comprise an unbroken run of integers from *1-n* up to *n*-1, and no others. Figure 19 shows an example of a minimal *magic* s-d square of order 6, the smallest yet found with this property. An interesting property of minimal s-d squares is that their bottom-right number (k) is even and occurs k/2 times in the square. The proof of this fact is due to Al Zimmermann. Let the sum of values in the last column and last row be C and R, respectively. In a minimal s-d square, the sum of the n-1 border values is 0, since all values between 1-n and n+1 occur exactly once. This sum is also C+R-k (k is counted twice), so C+R = k. Since C and R represent the number of occurrences of k, we also have C = R and so it must be that C = R = k/2. Hence k is even and occurs k/2 number of times.

|    |    |    |    |    |    |   |
|----|----|----|----|----|----|---|
| -5 | 0  | 3  | 3  | 3  | -3 | 1 |
| 1  | 1  | 1  | 3  | -1 | -4 | 1 |
| -1 | 2  | 1  | -1 | 2  | 2  | 5 |
| 1  | 4  | 1  | -1 | 0  | 0  | 5 |
| 0  | 0  | -2 | 2  | 2  | 3  | 5 |
| 5  | -5 | -2 | -1 | 1  | 4  | 2 |
| 1  | 1  | 2  | 2  | 5  | 7  | 2  2 |

Figure 19.
A minimal magic s-d square of order 6.

A further exotic case, again suggested by analogy with magic squares, is the *concentric* s-d square. In the example of Figure 20 both the 5×5 square as well as its inner 3×3 core form self-descriptive squares.

|    |    |    |    |    |   |
|----|----|----|----|----|---|
| 0  | 4  | 5  | -2 | -5 | 2 |
| 0  | 0  | 2  | 0  | -1 | 1 |
| 4  | -2 | -2 | 5  | 0  | 5 |
| -5 | 3  | 2  | -2 | 5  | 3 |
| 3  | -4 | -2 | 2  | 4  | 3 |
| 2  | 1  | 5  | 3  | 3  |   |

Figure 20. A concentric s-d square.

Computer searches have brought to light many examples of exotic s-d squares for orders *n*, $1 \leq n \leq 16$. In several cases we have been able to rule out the existence of various types of squares for certain orders. Our findings are summarized in Table 1.

| **Type / Order**       | 1 | 2 | 3 | 4 | 5 | 6 | 7 | 8 | 9 | 10 | 11 | 12 | 13 | 14 | 15 | 16 |
|------------------------|---|---|---|---|---|---|---|---|---|----|----|----|----|----|----|----|
| s-d                    | ✓ | ✗ | ✓ | ✓ | ✓ | ✓ | ✓ | ✓ | ✓ | ✓  | ✓  | ✓  | ✓  | ✓  | ✓  | ✓  |
| minimal s-d            | ✗ | ✗ | ✗ | ✗ | ✓ | ✓ | ✓ | ✓ | ✓ | ✓  | ✓  | ✓  | ✓  | ✓  | ✓  | ✓  |
| magic s-d              | ✓ | ✗ | ✗ | ✓ | ✓ | ✓ | ✓ | ✓ | ✓ | ✓  | ✓  | ✓  | ✓  | ✓  | ✓  | ✓  |
| minimal magic s-d      | ✗ | ✗ | ✗ | ✗ | ✗ | ✓ | ✓ | ✓ | ✓ | ✓  | ✓  | ✓  | ✓  | ✓  | ✓  | ✓  |
| perfect s-d            | ✓ | ✗ | ✗ | ✗ | ✓ | ✓ | ✓ | ✓ | ✓ | ?  | ?  | ?  | ?  | ?  | ?  | ?  |
| generic s-d            | ✗ | ✗ | ✗ | ✓ | ✓ | ✓ | ✓ | ✓ | ✓ | ?  | ?  | ?  | ?  | ?  | ?  | ?  |
| magic generic s-d      | ✗ | ✗ | ✗ | ✗ | ✓ | ✓ | ✓ | ? | ? | ?  | ?  | ?  | ?  | ?  | ?  | ?  |
| 2-variable generic s-d | ✗ | ✗ | ✗ | ✗ | ✓ | ✓ | ? | ? | ? | ?  | ?  | ?  | ?  | ?  | ?  | ?  |

Table 1. Summary of our findings. A tick indicates that an example has been found, a cross that no such square exists. Question marks identify cases where we don't know.

So much for a brief account of our investigations into self-descriptive squares. We cannot claim to have done more than identified some elementary features of the field and to have presented some interesting specimens. Two further topics that may become the subject of future interest are self-descriptive *cubes*, which are self-explanatory, and *co*-descriptive squares. Where a *self*-descriptive square can be viewed as a closed loop of length 1, so co-descriptive squares correspond to longer chains. In a co-descriptive loop of length 2 for example, we find a pair of squares, the row and column sums in each of which record the frequencies of their final elements occurring in the *other* square. In short, the squares describe eachother; see [3]. Meanwhile, should this topic prove of interest to any readers, we shall be glad to learn about any further advances made.

*Acknowledgements*

We are very grateful to Al Zimmermann for suggesting great improvements to our computer search algorithm and proving some useful properties about minimal s-d squares.

\*    \*    \*    \*    \*